\documentclass[11pt]{amsart}

\usepackage[english]{babel}
\usepackage{enumitem} 

\usepackage{graphicx}%
\usepackage{multirow}%
\usepackage[foot]{amsaddr}
\usepackage{amsmath,amssymb,amsfonts,setspace}%
\usepackage{bbold}
\usepackage{cite} 
\usepackage[title]{appendix}%
\numberwithin{equation}{section}
\usepackage{xcolor}%
\usepackage{textcomp}%
\usepackage{manyfoot}%
\usepackage{booktabs}%
\usepackage{algorithm}%
\usepackage{algorithmicx}%
\usepackage{algpseudocode}%
\usepackage{listings}%

\usepackage[colorlinks=true, linkcolor=blue, citecolor=red, urlcolor=blue]{hyperref}

\theoremstyle{thmstyleone}%
\newtheorem{theorem}{Theorem}%

\newtheorem{proposition}[theorem]{Proposition}%

\theoremstyle{thmstyletwo}%
\newtheorem{example}{Example}%
\newtheorem{remark}{Remark}%
\newtheorem{corollary}{Corollary}%
\newtheorem{lemma}{Lemma}%

\theoremstyle{thmstylethree}%
\newtheorem{definition}{Definition}%

\title{On Semi-Order Summability in Riesz Spaces}

\author{Mohammed El Amrani$^{1,*}$}
\address{$^1$Department of Mathematics, Faculty of Science, Laboratory of LAGA, Ibn Tofail University, B.P 242 Kenitra, Kenitra, 14000, Rabat-Sale-Kenitra, Morocco}
\email{mohammed.elamrani2@uit.ac.ma}
\thanks{$^*$Corresponding author.}

\author{Bouazza El Wahbi$^1$}
\email{bouazza.elwahbi@uit.ac.ma}

\author{Larbi Zraoula$^2$}
\address{$^2$Department of Mathematics, EREAM Team, LaREAMI-Lab, CRMEF-Kenitra, "Moroccan Association of Order Structures, Operators Theory, Application and sustainable Development (MAOSOTA)", B.P 242 Kenitra, Kenitra, 14000, Rabat-Sale-Kenitra, Morocco}
\email{zraoularbi@yahoo.fr}

\date{}

\begin{document}

\maketitle

\begin{abstract}
This paper examines the relationships between order summability, semi-order summability and norm summability in Riesz spaces and Banach lattices. First, we introduce the notion of semi-order summability, and establish its fundamental properties. We then prove that semi-order summability coincides with order summability in Dedekind complete Riesz spaces for arbitrary families, and in Dedekind $\sigma$-complete Riesz spaces for sequences. We also present new characterizations of laterally and $\sigma$-laterally complete Riesz spaces in terms of summability.

Furthermore, we examine the relationship between semi-order summability, and norm summability. We demonstrate that every semi-order-summable family is norm summable in normed Riesz spaces with order-continuous norm, and we derive new characterizations of Banach lattices with order-continuous norm. Finally, we prove that in a KB-space, a disjoint family is norm-summable if and only if its positive and negative parts are norm-summable, or equivalently if the family of its moduli is norm summable.
\end{abstract}

\maketitle

\section{Introduction}
The theory of summability has played a fundamental role in functional analysis since the early study of series in normed and Banach spaces. In Banach spaces, absolute and unconditional summability are among the most important concepts due to their deep connections with the topology of the underlying space, operator theory, and vector valued integration.
For a detailed exposition of summable families in normed vector spaces, the reader is referred to \cite{cd2, cd3, cd4, cd6}. In \cite{cd11}, the authors studied the norm summability of sequences in Banach lattices, but did not investigate the corresponding notion of order summability.

 Within the context of Riesz spaces and Banach lattices, where the order structure is important, it is logical to examine summability through order-theoretic concepts rather than solely through the norm. 
 The order summability of sequences (or, equivalently, order-convergent series) in a $\sigma$-Dedekind complete Riesz space was studied in \cite{cd7}.

The notion of order summability of families in Riesz spaces was introduced as an extension of the classical concept of summability from sequences and series to arbitrary families. It was subsequently generalized in \cite{cd5} to families indexed by an arbitrary nonempty set, whether countable or uncountable.

In this paper, we introduce a weaker notion of order summability, called \emph{semi-order summability} in Riesz space, We begin by establishing the basic properties of semi-order summable families. We then prove that the class of all semi-order summable families in a Riesz space $E$, indexed by a non-empty set $I$ forms a vector subspace of $E^I$ (Proposition~\ref{pp3}), each semi-order summable family is necessarily order bounded (Proposition~\ref{pp6}) and that every semi-order summable net whose directed index set has no maximal element is necessarily order convergent to zero (Proposition~\ref{p6}). We then establish a necessary and sufficient condition for a family in a Dedekind complete Riesz space to be semi-order summable, thereby obtaining a Cauchy-type criterion for semi-order summability (Theorem~\ref{t5}). 

Part of this paper is devoted to investigating the relationship between semi-order summability and order summability. First, we prove that every order summable family is semi-order summable with the same sum (Proposition~\ref{p5}). We then demonstrate that the converse is true in Dedekind-complete Riesz spaces, thereby showing that the two concepts are equivalent in this context. Specifically, they coincide for arbitrary families in Dedekind-complete Riesz spaces (Theorem~\ref{t6}) and for sequences in Dedekind $\sigma$-complete Riesz spaces (Proposition~\ref{pp8}). We also examine the special case of positive families and prove corresponding equivalence results (Proposition~\ref{p4}).

The paper further investigates disjoint families. We examine how lattice operations interact with order summation by proving that, for disjoint order summable families, taking positive and negative parts commutes with order summation (Theorem~\ref{t1}). We also investigate the relationship between order boundedness and the summability of disjoint families (Proposition~\ref{l2}, Proposition~\ref{l3}, and Corollary~\ref{cr1}). As an application of these results, we obtain two new characterizations of laterally complete and $\sigma$-laterally complete Riesz spaces (Proposition~\ref{pp1} and Proposition~\ref{pp2}).

The second part of the paper examines the relationship between semi-order summability and norm summability in normed Riesz spaces and Banach lattices. First, we prove that every norm summable net indexed by a directed set without maximal elements converges to zero in norm (Proposition~\ref{p2}). This extends a well-known property of convergent series to arbitrary families. We also show that, in a normed Riesz space, whenever a family is both norm summable and semi-order summable, the two notions produce the same sum (Proposition~\ref{p8}).

We then examine how the order continuity of the norm affects these summability notions.  We prove that every semi-order summable family is norm summable whenever the underlying normed Riesz space has an order-continuous norm (Proposition~\ref{p9}). As consequences, we obtain new characterizations of Banach lattices with order-continuous norms through the interplay of norm, order and semi-order summability of families (Theorem~\ref{t2}).

Finally, we study norm summability in KB-spaces. Although the norm summability of the positive and negative parts of a family always implies the norm summability of the family itself, the converse generally fails. We prove that the norm summability of a disjoint family is equivalent to the norm summability of its positive and negative parts, or equivalently, to the norm summability of the family of moduli (Theorem~\ref{tkb}).

For terminology and basic results on Riesz spaces and Banach lattices that are not recalled here, we refer the reader to  \cite{cd11, cd1, cd8}. Since order convergence is one of the fundamental notions underlying the results of this paper, the reader is referred to \cite{cd9, cd10} for a detailed account of its theory and basic properties.
\section{Summability of Families in Riesz Spaces}
In this section, we introduce the notion of semi-order summability for families in Riesz spaces and investigate its relationship with order summability. We also establish new characterizations of laterally complete and $\sigma$-laterally complete Riesz spaces in terms of order and semi-order summability.

\begin{definition} \label{d3}
Let $E$ be a Riesz space and let $I$ be a non-empty set. A family ${\left( x_{i} \right)}_{i \in I}$ of vectors in $E$ is said to be semi-order summable with sum $s$. If there exists a net $(y_\alpha)_{\alpha \in \Lambda}$ in $E^{+}$ such that $y_\alpha \downarrow 0$, and for each 
$\alpha \in \Lambda$ there exists a finite subset $J_\alpha \subseteq I$ such that, for every finite subset 
$J \subseteq I$,
$$J_\alpha \subseteq J \implies \left| \displaystyle \sum_{i \in J} x_i - s  \right| \le y_\alpha.$$
The vector $s$ is called the semi-order-sum of the family ${\left( x_{i} \right)}_{i \in I}$. In this case we write
$$\displaystyle \sum_{i \in I} {x_i}_{ / so} = s.$$
\end{definition} 
\textbf{Convention:} If $I = \varnothing$, we agree that
$\displaystyle \sum_{i \in I} {x_i}_{ / so} =  0.$

\begin{remark}
Let
$\Gamma=\biggl\{ \displaystyle \beta \mid \beta \subseteq I \text{ with } \beta \text{ is finite } \biggl\},$
ordered by inclusion, and define the net of finite partial sums by
$$X_\beta=\sum_{i\in \beta}x_i,\qquad \beta\in\Gamma.$$
Then the family $(x_i)_{i\in I}$ is semi-order summable with sum $s$ if and only if the net $(X_\beta)_{\beta\in\Gamma}$ converges in order to $s$. Thus, semi-order summability may be viewed as the order convergence of the net of finite partial sums.
\end{remark}

The next proposition shows that the collection of all semi-order summable families indexed by a nonempty set $I$ is closed under addition and scalar multiplication. Consequently, it forms a vector subspace of the product space $E^I$ .
\begin{proposition} \label{pp3}
Let $E$ be a Riesz space, let $I$ be a nonempty set, and let
$(x_i)_{i\in I}$ and $(y_i)_{i\in I}$ be two semi-order summable families in $E$.
 Then:
\begin{enumerate}[label=(\arabic*)]
    \item If $x_i \leq y_i $ for each $i \in I $, then
      $$\displaystyle \sum_{i \in I} {x_i}_{ / so} \leq \displaystyle \sum_{i \in I} {y_i}_{ / so}.$$
    
    \item $(x_i + y_i)_{i \in I} $ is a semi-order summable family, and
    $$\displaystyle \sum_{i \in I} {(x_i + y_i)}_{ / so} =
     \displaystyle \sum_{i \in I} {x_i}_{ / so} + \displaystyle \sum_{i \in I} {y_i}_{ / so}.$$
    
    \item $(\lambda x_i)_{i \in I}$ is a semi-order summable family, and
    $$\displaystyle \sum_{i \in I} {\lambda x_i}_{ / so} = 
    \lambda \displaystyle \sum_{i \in I} {x_i}_{ / so} \quad \text{for each} \ \lambda \in \mathbb{R}.$$
\end{enumerate}
\end{proposition}

The following proposition demonstrates that semi-order summability places strong restrictions on a family's order. Specifically, every semi-order summable family is necessarily order bounded.

\begin{proposition} \label{pp6}
Let $E$ be a Riesz space.
Every semi-order summable family in $E$ is order bounded.
\end{proposition}

\begin{proof}
Let $(x_\alpha)_{\alpha \in I}$ be a semi-order summable family in $E$ with sum $s$. 
So, there exists a net 
$(y_\beta)_{\beta \in \Lambda}$ in $E^+$ such that $y_\beta \downarrow 0$ and for each 
$\beta \in \Lambda$ there exists a finite subset 
$J_\beta \subseteq I$ 
such that, for every finite subset 
$J_{\beta} \subseteq J$, 
$$ \left| \displaystyle \sum_{i \in J} x_i - s  \right| \le y_{\beta}.$$
Fix $\beta_{0} \in \Lambda$. Then there exists a finite subset 
$J_{\beta_{0}} \subseteq I$ such that, 
$$ \left| \displaystyle \sum_{i \in J_{\beta_{0}}} x_i - s  \right| \le y_{\beta_{0}}.$$
Now let $k \in I  \setminus J_{\beta_{0}}$. Since $J_{\beta_{0}} \subseteq  J_{\beta_{0}} \cup \{ k \} $, it follows that,
$$\left| \displaystyle \sum_{i \in J_{\beta_{0}} \cup \{ k \}} x_i - s  \right| \le y_{\beta_{0}}.$$
So
$$\left| \displaystyle x_{k} + \sum_{i \in J_{\beta_{0}} } x_i - s  \right| \le y_{\beta_{0}}.$$
Therefore 
$$ -y_{\beta_{0}} + s - \sum_{i \in J_{\beta_{0}} } x_i \leq x_{k} \leq 
y_{\beta_{0}} + s - \sum_{i \in J_{\beta_{0}} } x_i $$
\end{proof}

\begin{proposition} \label{p6}
Let $(x_\alpha)_{\alpha \in I}$ be a net in a Riesz space $E$.
Suppose that the directed index set $I$ has no maximal element.\\
 If 
${\left( x_{\alpha}\right)_{\alpha \in I}}$ 
is a semi-order summable family, then $x_{\alpha}\xrightarrow{o} 0$.
\end{proposition}

\begin{proof}
Let $(x_\alpha)_{\alpha \in I}$ be a semi-order summable family in $E$ with sum $s$. 
So, there exists a net 
$(y_\beta)_{\beta \in \Lambda}$ in $E^+$ such that $y_\beta \downarrow 0$ and for each 
$\beta \in \Lambda$ there exists a finite subset 
$J_\beta \subseteq I$ 
such that, for every finite subset 
$J \subseteq I$, 
$$J_\beta \subseteq J \implies \left| \displaystyle \sum_{i \in J} x_i - s  \right| \le y_\beta.$$
On the other hand, since $I$ has no maximal element, there exists 
$\alpha_0 \in I$ such that $  \alpha_0  > i \quad \text{for every } i \in J_{\beta}$.
 So, for every 
$\alpha \ge \alpha_0$, we have 
\begin{align*}
\left| x_{\alpha}\right| &= \left| \displaystyle \sum_{i \in J_{\beta}\cup \{\alpha\}} x_i - \displaystyle \sum_{i \in J_{\beta}} x_i \right| \\
  &\leq \left| \displaystyle \sum_{i \in J_{\beta}\cup \{\alpha\}} x_i - s \right| + 
  \left| \displaystyle \sum_{i \in J_{\beta}} x_i - s \right| \\
  &\leq y_\beta+y_\beta= 2y_\beta
\end{align*}
Therefore, $x_{\alpha}\xrightarrow{o} 0$.
\end{proof}

The assumption that the directed index set has no maximal element in the preceding proposition is essential. In general, a semi-order summable net need not converge to zero in order. The following example illustrates this fact.
\begin{example}
Let $E=\mathbb{R}$ and let $I=]0,1]$, directed by the usual order. Define the net $(x_\alpha)_{\alpha\in I}$ by
$$
x_\alpha=
\begin{cases}
\alpha^{2}, & \text{if } \alpha=\dfrac{1}{n}\text{ for some }n\in\mathbb{N}^{*},\\[4pt]
0, & \text{otherwise}.
\end{cases}
$$
Then $(x_\alpha)_{\alpha\in I}$ is a semi-order summable family in $E$, and
$$\sum_{\alpha\in I}x_\alpha
=\sum_{n=1}^{\infty}\frac{1}{n^{2}}
=\frac{\pi^{2}}{6}.$$
However, $(x_\alpha)_{\alpha\in I}$ does not converge to $0$ in order. Indeed,
$$x_\alpha\xrightarrow{o}1,$$
since $1$ is the maximal element of the directed set $I$.
\end{example}
\begin{corollary}
Let $(x_\alpha)_{\alpha \in I}$ be a semi-order summable net in a Riesz space $E$.
Then $x_{\alpha}\xrightarrow{o} 0$  if and only if one of the following conditions is satisfied:
\begin{enumerate}[label=(\arabic*)]
\item The index set $I$ has no maximal element.
\item The index set $I$ has a maximal element $\alpha_0$ and $x_{\alpha_0}=0$.
\end{enumerate}
\end{corollary}

\begin{proposition} \label{pp9}
Let $E$ be a Riesz space and $(x_i)_{i \in I}$ be a semi-order summable family in $E$. 
Then there exists a net $(z_\alpha)_{\alpha \in \Lambda}$ in $E^{+}$ such that 
$z_\alpha \downarrow 0$ and, for each $\alpha \in \Lambda$, there exists a finite subset 
$J_\alpha \subseteq I$ such that
$$\left| \sum_{i \in J} x_i  \right| \leq z_\alpha$$
for all finite subsets $J \subseteq I$ satisfying $J \cap J_\alpha = \varnothing$.
\end{proposition}

\begin{proof}
Let $(x_i)_{i \in I}$ be a semi-order summable family in $E$ with sum $s$. 
Then there exists a net $(y_\alpha)_{\alpha \in \Lambda}$ in $E^{+}$ such that 
$y_\alpha \downarrow 0$ and, for each $\alpha \in \Lambda$, there exists a finite subset 
$J_\alpha \subseteq I$ such that
$$\left| \sum_{i \in J_{\alpha}} x_i - s \right| \le y_\alpha$$
Now let $J \subseteq I$ be a finite subset such that $J \cap J_\alpha = \varnothing$. 
Since $J_{\alpha} \subseteq  J \cup J_{\alpha}$, it follows that
\begin{align*}
\left| \sum_{i \in J} x_i  \right| &= \left| \sum_{i \in J \cup J_{\alpha}} x_i - \sum_{i \in J_{\alpha}} x_i \right| \\
  &\leq \left| \sum_{i \in J \cup J_{\alpha}} x_i -s  \right| + 
  \left| \sum_{i \in J_{\alpha}} x_i -s  \right| \\
  &\leq y_\alpha + y_\alpha= 2y_\alpha
\end{align*}
Hence, it is enough to choose $z_{\alpha} = 2y_{\alpha}$ for each $\alpha \in \Lambda$.
\end{proof}

The following theorem provides a Cauchy-type characterization of semi-order summable families in Dedekind complete Riesz spaces. More precisely, it demonstrates that the semi-order summability of a family is entirely dependent on the order behavior of the tails of its finite partial sums.
\begin{theorem} \label{t5}
Let $E$ be a Dedekind complete Riesz space, and let $(x_i)_{i \in I}$ be a family in $E$. Then the following assertions are equivalent:
\begin{enumerate}[label=(\arabic*)]
    \item $(x_i)_{i \in I}$ is semi-order summable in $E$.
    \item There exists a net $(z_\alpha)_{\alpha \in \Lambda}$ in $E^{+}$ such that 
$z_\alpha \downarrow 0$ and, for each $\alpha \in \Lambda$, there exists a finite subset 
$J_\alpha \subseteq I$ such that
$$\left| \sum_{i \in J} x_i  \right| \leq z_\alpha$$
for all finite subsets $J \subseteq I$ satisfying $J \cap J_\alpha = \varnothing$.
\end{enumerate}
\end{theorem}

\begin{proof}
$(1)\Rightarrow (2)$ By Proposition~\ref{pp9}. \\
$(2)\Rightarrow (1)$ Let $\Gamma=\biggl\{ \displaystyle \beta \mid \beta \subseteq I \text{ with } \beta \text{ is finite } \biggl\}$ and 
$X_{\beta}=\displaystyle \sum_{i\in \beta} x_{i}$ for each $\beta \in \Gamma$.
Assume that $(2)$ holds, then there exists a net $(z_\alpha)_{\alpha \in \Lambda}$ in $E^{+}$ such that 
$z_\alpha \downarrow 0$ and, for each $\alpha \in \Lambda$, there exists a finite subset 
$J_\alpha \subseteq I$ such that
$$\left| \sum_{i \in J} x_i  \right| \leq z_\alpha$$
for all finite subset $J \subseteq I$ satisfying $J \cap J_\alpha = \varnothing$.\\
Now let $\beta, \gamma \in \Gamma$ be such that $J_\alpha \subseteq \beta$ and $J_\alpha \subseteq \gamma$. 
Then $\left( \beta \setminus \gamma \right) \cap J_{\alpha} = \varnothing$ and 
$\left( \gamma \setminus \beta \right) \cap J_{\alpha} = \varnothing$, so
\begin{align*}
\left| X_{\beta} - X_{\gamma}  \right| &= \left| \sum_{i \in \beta } x_i - \sum_{i \in \gamma} x_i \right| \\
  &= \left| \sum_{i \in \beta \setminus \gamma } x_i - \sum_{i \in \gamma \setminus \beta} x_i \right| \\
  & \leq \left|  \sum_{i \in \beta \setminus \gamma } x_i  \right| + 
  \left| \sum_{i \in \gamma \setminus \beta} x_i  \right| \\
  & \leq z_\alpha + z_\alpha = 2z_\alpha.
\end{align*}
Hence, ${(X_{\beta})}_{\beta \in \Gamma}$ is order Cauchy. Since $E$ is Dedekind complete, it follows from \cite[Proposition 2.3]{cd14} that ${(X_{\beta})}_{\beta \in \Gamma}$ is order convergent. Therefore, the family $(x_i)_{i \in I}$ is semi-order summable.
\end{proof}

\begin{definition} \cite{cd5} \label{d2}
Let $E$ be a Riesz space and let $I$ be a non-empty set. A family ${\left( x_{i} \right)}_{i \in I}$ of vectors in $E$ is said to be order summable with sum $s$ 
if both subsets
$$A^{+} = \left\{\displaystyle \sum_{i \in J} {x_{i}}^{+} \mid J \subseteq I \text{ with } J \text{ is finite } \right\}$$
and
$$A^{-} = \left\{\displaystyle \sum_{i \in J} {x_{i}}^{-} \mid J \subseteq I \text{ with } J \text{ is finite } \right\}$$
have suprema in $E^{+}$, denoted by $s_{1}$ and $s_{2}$ respectively, and $s=s_{1}-s_{2}$. 
The vector $s$ is called the order-sum of the family ${\left( x_{i} \right)}_{i \in I}$. In this case we write
$$\displaystyle \sum_{i \in I} {x_i}_{ / o} = s.$$
\end{definition}
\textbf{Convention:} If $I = \varnothing$, we agree that
$\displaystyle \sum_{i \in I} {x_i}_{ / o} =  0.$

\begin{proposition} \cite{cd5} \label{p3}
Let \( E \) be a Dedekind complete Riesz space, and let \( x = (x_i)_{i \in I} \subseteq E \) and \( y = (y_i)_{i \in I} \subseteq E \) be two order summable families. Then:

\begin{enumerate}[label=(\arabic*)]
    \item If $x_i \leq y_i $ for each $i \in I $, then
      $$\displaystyle \sum_{i \in I} {x_i}_{ / o} \leq \displaystyle \sum_{i \in I} {y_i}_{ / o}$$.
    
    \item $(x_i + y_i)_{i \in I} $ is an order summable family, and
    $$\displaystyle \sum_{i \in I} {(x_i + y_i)}_{ / o} =
     \displaystyle \sum_{i \in I} {x_i}_{ / o} + \displaystyle \sum_{i \in I} {y_i}_{ / o}.$$
    
    \item $(\lambda x_i)_{i \in I}$ is a order summable family, and
    $$\displaystyle \sum_{i \in I} {\lambda x_i}_{ / o} = 
    \lambda \displaystyle \sum_{i \in I} {x_i}_{ / o} \quad \text{for each} \ \lambda \in \mathbb{R}.$$
\end{enumerate}
\end{proposition}

The following proposition establishes the equivalence between order summability and semi-order summability for positive families.
\begin{proposition} \label{p4}
Let $E$ be a Riesz space. A positive family $(x_i)_{i \in I}$ in $E$
is semi-order summable if and only if it is order summable.\\
Moreover, in the case of summability, we have
$$ \displaystyle \sum_{i \in I} {x_i}_{ / so} = \sum_{i \in I} {x_i}_{ / o}.$$
\end{proposition}

\begin{proof}
Let ${\left( x_{i} \right)}_{i \in I}$ be a family in $E^{+}$, 
$\Lambda=\biggl\{ \displaystyle \alpha \mid \alpha \subseteq I \text{ with } \alpha \text{ is finite } \biggl\}$ and 
$z_{\alpha}=\displaystyle \sum_{i\in \alpha} x_{i}$ for each $\alpha \in \Lambda$.
It now suffices to observe that the family $(x_i)_{i \in I}$ is order summable with sum $s$ if and only if the net $z_\alpha \uparrow s$.
\end{proof}

The following proposition demonstrates that every order summable family is also semi-order-summable, and that the corresponding sums are identical.
\begin{proposition} \label{p5}
Let $E$ be a Riesz space and let $(x_i)_{i \in I}$ be an order summable family in $E$. Then $(x_i)_{i \in I}$ is semi-order summable in $E$. Moreover, we have 
$$ \displaystyle \sum_{i \in I} {x_i}_{ / so} = \sum_{i \in I} {x_i}_{ / o}.$$
\end{proposition}

\begin{proof}
Let $(x_i)_{i \in I}$ be an order summable family in $E$ with sum $s$. Then the families
$({x_i}^{+})_{i \in I}$ and $({x_i}^{-})_{i \in I}$ are order summable,
with respective sums $s_1$ and $s_2$ such that $s=s_{1}-s_{2}$. 
So, by Proposition~\ref{p4}, $({x_i}^{+})_{i \in I}$ and $({x_i}^{-})_{i \in I}$ are semi-order summable, with respective sums $s_1$ and $s_2$. It now suffices to apply Proposition~\ref{pp3}.
\end{proof}

As shown by the following example, semi-order summability does not generally imply order summability.

\begin{example} \label{ex1} 
Let $c_{00} = \{ (x_n)_{n \in \mathbb{N}} \subseteq {\mathbb{R}}^{\mathbb{N}} \mid x_n = 0 \text{ for all but finitely many } n \}$, $v=(1, 2, 3, ...) \in {\mathbb{R}}^{\mathbb{N}}$, 
$(e_n)_{n \in \mathbb{N^{*}}} \subseteq c_{00}$ denotes the family of unit vectors, 
$\mathrm{span}\{v\} = \{ \lambda v \mid \lambda \in \mathbb{R} \}$,
$E=c_{00}\bigoplus \mathrm{span}\{v\}$, and
$(x_n)_{n \in \mathbb{N^{*}}} \subseteq E$ such that $x_{n}=e_{n}-e_{n+1}$ for each $n \in N^{*}$.\\
$(x_n)_{n \in \mathbb{N^{*}}}$ is a semi-order summable family in $E$. Indeed, let
$(y_n)_{n \in \mathbb{N^{*}}} \subseteq E$ such that 
$y_{n}=v- \displaystyle \sum_{k =1}^{n} {k.e_{k}}$ for each $n \in N^{*}$, it is clear that
$y_n \downarrow 0$ in $E$.\\
Let $n \in \mathbb{N^{*}}$ and define 
$J_n = \{1, 2, \dots, n\}$. 
Let $J \subseteq \mathbb{N}$ be finite. 
If $J_n \subseteq J$ , then we can write
$J = \{1, 2, \dots, n, n_1, n_2, n_3, \dots\},$
where $n_1, n_2, \dots$ denote the remaining elements of $J$ beyond $n$. So,
\begin{align*}
\left|  \displaystyle \sum_{i \in J} {x_{i}} - e_{1} \right| &= \left| -e_{2} + e_{2}- e_{3}+ e_{3} + \dots + e_{n}-e_{n+1}+ e_{n_{1}}- e_{n_{1} +1} \dots \right| \\
  &=\left|  -e_{n+1}+ e_{n_{1}}- e_{n_{1} +1} + e_{n_{2}} - e_{n_{2} +1} + \dots  \right| \\
  &= (0, 0, \dots, \underset{\text{n-th}}{0}, a_{n_1}, a_{n_2}, \dots), \text{where each} \,
 a_{n_i} \in \{0, 1 \}\\
 & \leq y_{n}.
\end{align*}
However, $(x_n)_{n \in \mathbb{N^{*}}}$ is not order summable, since 
$$A^{+} = \left\{\displaystyle \sum_{n \in J} {x_{n}}^{+} \mid J \subseteq  \mathbb{N^{*}} \text{ with } J \text{ is finite } \right\}$$
and
$$A^{-} = \left\{\displaystyle \sum_{n \in J} {x_{n}}^{-} \mid J \subseteq  \mathbb{N^{*}} \text{ with } J \text{ is finite } \right\}$$
 do not have a supremum in $E$.
\end{example}

\begin{example}
Consider $\mathbb{R}^{2}$ equipped with the standard order. 
Let $I=\{1, 2\}$, $x_{1}=(3, -2)$, and $x_{2}=(-1, 4)$. 
The power set of $I$ is $\{\emptyset, \{1\}, \{2\}, \{1, 2\}\}$.\\
We have 
$$\bigvee_{K \subseteq \{1,2\}} \left\{ \sum_{i \in K} x_i \right\} = (0,0) \vee (3, -2) \vee (-1, 4) \vee (2, 2) = (3, 4)$$ 
On the other hand 
$$\sum_{i \in I} x_{i}^+ = x_1^+ + x_2^+ = (3, 0) + (0, 4) = (3, 4)$$Therefore:$$\sum_{i \in I} x_{i}^+ = \bigvee_{K \subseteq \{1,2\}} \left\{ \sum_{i \in K} x_i \right\}$$
This observation is generalized by the following lemma.
\end{example}

\begin{lemma} \label{l4}
Let $E$ be a Riesz space, $I$ a finite index set, and $(x_i)_{i \in I} \subseteq E$ a finite family. Then 
$$\displaystyle \sum_{i \in I} {x_{i}}^{+}= \bigvee_{K \subseteq I} \left\{ \sum_{i \in K} {x_{i}} \right\}$$
\end{lemma}

\begin{proof}
If $I = \emptyset$, then the result is trivial.\\
Now without loss of generality, we may assume that $I = \{1,2,\dots,n\}$. Then let us prove that 

$$\displaystyle \sum_{i=1}^{n} x_i^{+} = 
\bigvee_{K \subseteq \{1,2,\dots,n\}} \left\{ \sum_{i \in K} x_i \right\}$$ 
by induction on $n$.\\
For $n=1$, the collection of subsets of $I=\{1\}$ is 
$\{\varnothing, \{1\}\}.$ The identity becomes ${x_{1}}^{+} = \bigvee \{0, x_{1}\}.$\\
For $n=2$, we have
\begin{align*}
{{x}_{1}}^{+} + {{x}_{2}}^{+}&= x_{1}\vee 0 + x_{2}\vee 0 \\
  &=  \left( \left( x_{1}\vee 0 \right) + x_{2} \right) \vee \left( x_{1}\vee 0 \right) \\
  &= \left( x_{1} + x_{2} \right) \vee x_{2} \vee x_{1} \vee 0
\end{align*} 
Since
$$ \displaystyle \left( x_{1} + x_{2} \right) \vee x_{1} \vee x_{2} \vee 0 =
\bigvee_{K \subseteq \{1,2\}} \left\{ \sum_{i \in K} x_i \right\},$$
the identity is verified for $n=2$.\\
Assume that the formula holds for $n$ elements. We prove it for $n+1$.\\
Let
$\displaystyle S_{n} = \sum_{i=1}^{n} {x_i}^{+}.$
By the induction hypothesis,
$\displaystyle S_{n} = \bigvee_{K \subseteq \{1,\dots,n\}}
\left\{ \sum_{i \in K} x_i \right\}.$\\ 
We have $S_{n+1}=S_{n}+{x_{n+1}}^{+}$. Let us set $x_{1} = S_{n}$ and $x_{2} = x_{n+1}$, using the result for the case $n=2$, we have
$S_{n+1} = \left( S_{n} + x_{n+1} \right) \vee x_{n+1} \vee S_{n} \vee 0.$\\
Since $S_{n} \geq 0$ and $S_{n} + x_{n+1}  \geq x_{n+1}$, it follows that
\begin{align*}
S_{n+1}&= \left( S_{n} + x_{n+1} \right) \vee s_{n} \\
&= \left( \bigvee_{K \subseteq \{1,\dots,n\}} \left\{ \sum_{i \in K} x_i \right\} + x_{n+1} \right) \vee \left( \bigvee_{K \subseteq \{1,\dots,n\}} \left\{ \sum_{i \in K} x_i \right\} \right) \\
  &= \left( \bigvee_{K \subseteq \{1,\dots,n\}} \left\{ \sum_{i \in K} x_i  + x_{n+1} \right\}  \right) \vee \left( \bigvee_{K \subseteq \{1,\dots,n\}} \left\{ \sum_{i \in K} x_i \right\} \right)   \\
  &= \bigvee_{K \subseteq \{1,\dots, n+1\}} \left\{ \sum_{i \in K} x_i \right\}
\end{align*} 
Hence, the proof is complete.
\end{proof}

The following theorem establishes that the notions of order summability and semi-order summability coincide in Dedekind complete Riesz spaces.
\begin{theorem} \label{t6}
Let $(x_i)_{i \in I}$ be a family in a Dedekind complete Riesz space $E$. 
Then $(x_i)_{i \in I}$ is order summable if and only if it is semi-order summable.
\end{theorem}

\begin{proof}
$ \Rightarrow )$ The result follows from Proposition~\ref{p5}.\\
$ \Leftarrow )$ 
Let $(x_{i})_{i \in I}$ be a semi-order summable family in $E$ with sum $s$. 
So, there exists a net 
$(y_\alpha)_{\alpha \in \Lambda}$ in $E^+$ such that $y_\alpha \downarrow 0$ and for each 
$\alpha \in \Lambda$ there exists a finite subset 
$J_\alpha \subseteq I$ 
such that, for every finite subset 
$J_{\alpha} \subseteq J$, 
$$ \left| \displaystyle \sum_{i \in J} x_i - s  \right| \le y_{\alpha}.$$
Fix $\alpha_{0} \in \Lambda$. Then there exists a finite subset 
$J_{\alpha_{0}} \subseteq I$ such that, 
$$ \left| \displaystyle \sum_{i \in J_{\alpha_{0}}} x_i - s  \right| \le y_{\alpha_{0}}.$$
Now let $K$ be a finite subset of $I$. Then 
$$ \left| \displaystyle \sum_{i \in K \cup J_{\alpha_{0}}} x_i - s  \right| \leq y_{\alpha_{0}}.$$
So
$$ \sum_{i \in K} x_i + \sum_{i \in J_{\alpha_0} \setminus K} x_i \leq y_{\alpha_0} + s.$$
Hence 
\begin{align*}
\sum_{i \in K} x_i &\leq  y_{\alpha_0} + s - \sum_{i \in J_{\alpha_0} \setminus K} x_i \\
  &\leq  y_{\alpha_0} + s + \sum_{i \in J_{\alpha_0} \setminus K} | x_i | \\
  & \leq y_{\alpha_0} + s + \sum_{i \in J_{\alpha_0}} | x_i |. 
\end{align*}
Therefore, the family of finite partial sums
$ \displaystyle \left\{ \sum_{i \in K} x_i \;|\; K \subseteq I \text{ finite} \right\} $
is bounded above in $E$. Then, by Lemma~\ref{l4}, the set
$ \displaystyle \left\{ \sum_{i \in K} {x_i}^{+} \;|\; K \subseteq I \text{ finite} \right\} $
is also bounded above. Since $E$ is Dedekind complete, it admits a supremum.\\ 
To show that
$ \displaystyle \left\{ \sum_{i \in K} {x_i}^{-} \;|\; K \subseteq I \text{ finite} \right\} $
admits a supremum, it suffices to observe that the family $(-x_i)_{i \in I}$ is summable and that
${x}^{-} = {(-x)}^{+}.$\\
 Thus, the proof is finished.
\end{proof} 

The following proposition gives the sequential version of the preceding result. It shows that order summability and semi-order summability coincide in $\sigma$-Dedekind complete Riesz spaces.
\begin{proposition} \label{pp8}
Let $(x_n)_{n \in \mathbb{N}}$ be a sequence in a $\sigma$-Dedekind complete Riesz space $E$. 
Then $(x_n)_{n \in \mathbb{N}}$ is order summable if and only if it is semi-order summable.
\end{proposition}

\begin{proof}
The proof is similar to that of Theorem~\ref{t6}.
\end{proof}

The following theorem proves that the operations of taking the positive and negative parts of disjoint order summable families commute with order summation.
\begin{theorem} \label{t1}
Let $(x_i)_{i \in I}$ be a disjoint order summable family in a Riesz space $E$ with sum $s$, such that
 $\displaystyle \sum_{i \in I} {{x_i}^{+}}_{ / o}=s_{1}$ and 
 $\displaystyle \sum_{i \in I} {{x_i}^{-}}_{ / o}=s_{2}$. Then 
  \begin{center}
  $s^{+} = s_{1}$ and $s^{-} = s_{2}$.
  \end{center}
\end{theorem}

\begin{proof}
First, we note that if the family $(x_i)_{i \in I}$ is disjoint, then ${x_i}^{+} \wedge {x_j}^{-} = 0$ for all $i,j \in I$.\\ 
Let 
$\Gamma=\biggl\{ \displaystyle \beta \mid \beta \subseteq I \text{ with } \beta \text{ is finite } \biggl\}$
$X_{\beta}= \displaystyle \sum_{i \in \beta } {x_i}^{+}$, and  
$Y_{\beta}= \displaystyle \sum_{i \in \beta } {x_i}^{-}$.
Then, by \cite[Lemma 1.4]{cd1}, we have 
$$0 \leq X_{\beta} \wedge Y_{\beta} \leq \displaystyle \sum_{i,j \in \beta } {x_i}^{+} \wedge {x_j}^{-}=0 
\quad \text{for each } \beta \in \Gamma$$ 
Since ${(x_i)}_{i \in I}$ is order summable, it follows that
$$X_{\beta}\xrightarrow{o} s_{1} \quad \text{and} \quad Y_{\beta}\xrightarrow{o} s_{2} $$
So 
$$ X_{\beta} \wedge Y_{\beta} \xrightarrow{o} s_{1}\wedge s_{2}$$ 
Therefore $s_{1}\wedge s_{2}=0$; hence, $s^{+} = s_{1}$ and $s^{-} = s_{2}$.
\end{proof}

The disjointness assumption in the preceding theorem is essential. In general, the operations of taking positive and negative parts do not commute with order summation. The following example illustrates this fact.

\begin{example}
Let $E=\mathbb{R}^{2}$ endowed with the coordinatewise order, and let $I=\{1,2\}$. Define
$$x_{1}=(-2,1), \qquad x_{2}=(3,-5).$$
Then $(x_i)_{i\in I}$ is an order summable family with
$$\sum_{i\in I}x_i=(1,-4).$$
Moreover,
$$\left(\sum_{i\in I}x_i\right)^+=(1,0)\neq(3,1)=\sum_{i\in I}x_i^+,$$
and
$$\left(\sum_{i\in I}x_i\right)^-=(0,4)\neq(2,5)=\sum_{i\in I}x_i^-.$$
Hence, in general, taking positive and negative parts does not commute with order summation.
\end{example}

The next proposition establishes that in a Dedekind complete Riesz space, every order bounded disjoint family of positive elements is order summable, and that its order sum coincides with its supremum.
\begin{proposition} \label{l2}
Let $E$ be a Dedekind complete Riesz space, and let ${\left( x_{i}\right)_{i \in I}}$ be an order bounded disjoint family in $E^+$. Then
${\left( x_{i}\right)_{i \in I}}$ is order summable and
$$ \sum_{i\in I} {x_i}_{ / o} = \bigvee_{i\in I} x_i. $$
\end{proposition}

\begin{proof}
Since $E$ is Dedekind complete and the family $(x_i)_{\i \in I}$ is order bounded, then 
$\displaystyle \bigvee_{i\in I} x_i$
 exists in $E$.\\
 For any finite subset $J \subseteq I$, we have 
 $$\displaystyle \sum_{i\in J} x_i=\bigvee_{i\in J} x_i \leq \bigvee_{i\in I} x_i.$$
 This implies that $ \displaystyle \bigvee_{i\in I} x_i $ is an upper bound of 
 $A = \left\{\displaystyle \sum_{i \in J} x_{i} \mid J \subseteq I \text{ with } J \text{ is finite } \right\}$.\\
 If $s$ is an upper bound of $A$, then $s \geq x_{i}$ holds for each $i \in I$.\\
Consequently, $s \geq \displaystyle \bigvee_{i\in I} x_i $. 
Hence $s$ is the least upper bound of $A$.\\
This completes the proof.
\end{proof} 

The following proposition is the sequential counterpart of the preceding result. It shows that, in a $\sigma$-Dedekind complete Riesz space, every order bounded disjoint sequence of positive elements is order summable, and that its order sum coincides with its supremum.
\begin{proposition} \label{l3}
Let $E$ be a Dedekind $\sigma$-complete Riesz space, and let ${\left( x_{n}\right)_{n \in \mathbb{N}}}$ be an order bounded disjoint sequence in $E^+$. Then
${\left( x_{n}\right)_{n \in \mathbb{N}}}$ is order summable and
$$\displaystyle \sum_{n\in \mathbb{N}} {x_n}_{ / o} = \bigvee_{n\in \mathbb{N}} x_n. $$
\end{proposition}

\begin{proof}
The proof is similar to that of Proposition~\ref{l2}.
\end{proof}

\begin{corollary} \label{cr1}
Let $E$ be a Dedekind complete (respectively, $\sigma$-Dedekind complete) Riesz space, and let ${\left( x_{i}\right)_{i \in I}}$ be an order bounded disjoint family (respectively, order bounded sequence) in $E$. Then
${\left( x_{i}\right)_{i \in I}}$ is order summable.
\end{corollary}

\begin{proof}
It is enough to verify that ${\left({x_{i}}^{+}\right)}_{i \in I}$ and 
${\left({x_{i}}^{-}\right)}_{i \in I}$ are disjoint and to apply
Propositions~\ref{l2} and~\ref{l3}.
\end{proof}

The following proposition gives equivalent characterizations of laterally complete Riesz spaces in terms of the order summability and semi-order summability of disjoint families.
\begin{proposition} \label{pp1}
Let $E$ be a Riesz space. Then the following assertions are equivalent:
\begin{enumerate}[label=(\arabic*)]
    \item $E$ is laterally complete.
    \item Every disjoint family in $E$ is order summable.
    \item Every disjoint family in $E$ is semi-order summable.
\end{enumerate}
\end{proposition}

\begin{proof}
$\text{(1)} \Rightarrow \text{(2)}$ 
The proof of this implication is analogous to that of  Proposition~\ref{l2}.\\
$\text{(2)} \Rightarrow \text{(3)}$ This is trivial by Proposition~\ref{p5}.\\
$\text{(3)} \Rightarrow \text{(1)}$ 
Let $A$ be a subset of disjoint vectors in $E^{+}$ and let ${\left( x_{a}\right)_{a \in A}}$ be a family such that 
$x_{a}=a$ for each $a \in A$. The family ${\left( x_{a}\right)_{a \in A}}$
 is disjoint and therefore semi-order summable with sum $s$. Hence, by Proposition~\ref{p4}, 
 it is order summable with the same sum $s$.\\
 It is clear that $a \leq s$ for each $a \in A$. We now consider $r \in E$ to be an upper bound of $A$.
  So we have that 
   $$\sum_{a\in J} x_a=\bigvee_{a\in J} x_a \leq r \quad \text{for each finite subset} \quad J \subseteq A.$$
   It follows that $s\leq r$, hence $s$ is the least upper bound of $A$.
This completes the proof.
\end{proof}

The following proposition provides the sequential counterpart of the preceding characterization. It gives equivalent characterizations of $\sigma$-laterally complete Riesz spaces in terms of the order summability and semi-order summability of disjoint sequences.
\begin{proposition} \label{pp2}
Let $E$ be a Riesz space. Then the following assertions are equivalent:
\begin{enumerate}[label=(\arabic*)]
    \item $E$ is $\sigma$-laterally complete.
    \item Every disjoint sequence in $E$ is order summable.
    \item Every disjoint sequence in $E$ is semi-order summable.
\end{enumerate}
\end{proposition}

\begin{proof}
The proof is similar to that of Proposition~\ref{pp1}.
\end{proof}
\section{Summability of Families in Banach Lattices}

In this section, we examine the relationship between semi-order summability and norm summability in Riesz spaces and Banach lattices. We also provide new characterizations of Banach lattices with order-continuous norms. 

\begin{definition} \cite{cd6} \label{d1}
Let $X$ be a normed vector space. A family ${\left( x_{i} \right)}_{i \in I}$ of vectors in $X$ is said to be
norm summable if there exists a vector $s \in X$ such that for every $\varepsilon > 0$,
there exists a finite subset $J_\varepsilon \subseteq I$ with the property that, for every finite
subset $J \subseteq I$,
$$J_\varepsilon \subseteq J \implies \left\| \displaystyle \sum_{i \in J} x_i - s \right\| \le \varepsilon.$$
The vector $s$ is called the norm-sum of the family ${\left( x_{i} \right)}_{i \in I}$, and in this case we write
$$\displaystyle \sum_{i \in I} {x_i}_{ / \|.\|} = s.$$
\end{definition}

\textbf{Convention:} If $I = \varnothing$, we agree that
$\displaystyle \sum_{i \in I} {x_i}_{ / \|.\|} =  0.$

\begin{remark} \label{r1}
Let ${\left( x_{i} \right)}_{i \in I}$ be a family in a normed vector space $X$,\\ 
$\Gamma=\biggl\{ \displaystyle \beta \mid \beta \subseteq I \text{ with } \beta \text{ is finite } \biggl\}$ and 
$X_{\beta}=\displaystyle \sum_{i\in \beta} x_{i}$ for each $\beta \in \Gamma$, so 
${\left( x_{i} \right)}_{i \in I} $ is norm summable if and only if there exists $s \in X$ such that 
$X_\beta \xrightarrow{\|.\|} s$.
\end{remark}

\begin{remark}
Let $(x_i)_{i\in I}$ and $(y_i)_{i\in I}$ be two norm-summable families in a normed Riesz space. If
$x_i\le y_i$ for each $ i\in I,$
then
$$\sum_{i\in I}x_i{}_{/\|\cdot\|}
\le
\sum_{i\in I}y_i{}_{/\|\cdot\|}.$$
\end{remark}

\begin{remark} \label{r2}
If $(x_i)_{i \in I}$ is a norm summable family in a normed vector space $X$, then by \cite[Corollary~5.28]{cd6} the set $\{i \in I \mid x_i \neq 0\}$ is countable.
\end{remark}

\begin{proposition}\cite{cd6} \label{p1}
Let $X$ be a normed vector space, and let $x = (x_i)_{i \in I} \subseteq X$ and 
$y = (y_i)_{i \in I} \subseteq X$ be two norm summable families. Then: 

\begin{enumerate}[label=(\arabic*)]
    \item The family $(x_i + y_i)_{i \in I}$ is norm summable, and 
    $$\displaystyle \sum_{i \in I} {(x_i + y_i)}_{ / \|.\|}  =\displaystyle \sum_{i \in I} {x_i}_{ / \|.\|} + \displaystyle \sum_{i \in I} {y_i}_{ / \|.\|}.$$
    \item For each $\lambda \in \mathbb{R}$, the family $(\lambda x_i)_{i \in I}$ is norm summable, and
    $$ \displaystyle \sum_{i \in I} {(\lambda x_i)}_{ / \|.\|} = 
    \lambda \displaystyle \sum_{i \in I} {x_i}_{ / \|.\|}.$$
\end{enumerate}
\end{proposition}

\begin{proposition} \label{l5}
Let $X$ be a normed vector space, and let $(x_i)_{i\in I}$ be a norm summable family in $X$. Then there exists
 $M \in \mathbb{R}^+$ such that
$\displaystyle  \left\|\sum_{i\in J} x_i\right\| \le M$
for every finite subset $J \subseteq I$.
\end{proposition}

\begin{proof}
Let $(x_i)_{i\in I}$ be a norm summable family in $X$ such that
$\displaystyle \sum_{i\in I} {x_{i}}_{ / \|.\|} = s.$ 
Then there exists a finite subset $J_{1} \subseteq I$ such that
$\displaystyle \left\|\sum_{i\in J} x_i-s\right\|\leq 1 $
for every finite subset $J\subseteq I$ satisfying $J_{1}\subseteq J.$ 
Now let $K$ be a finite subset of $I$. Since
$J_{1}\subseteq K\cup J_{1},$
it follows that
$$ \displaystyle \left\| \sum_{i\in K \cup J_{1}} x_i  \right\| - \left\| s \right\| \leq
 \left| \left\| \sum_{i\in K \cup J_{1}} x_i  \right\| - \left\| s \right\| \right| \leq 
\left\| \sum_{i\in K \cup J_{1}} x_i - s  \right\| \leq 1 .$$
So, $$ \displaystyle \left\|  \sum_{i\in K \cup J_{1}} x_i \right\| \leq 1 + \left\| s \right\|.$$
On the other hand,
\begin{align*}
\left\| \sum_{i\in K } x_i \right\| - \left\| \sum_{i\in J_{1} \setminus K} x_i \right\| &\leq 
\left| \left\| \sum_{i\in K } x_i \right\| - \left\| \sum_{i\in J_{1} \setminus K} x_i \right\|  \right| \\
 &= \left| \left\| \sum_{i\in K } x_i \right\| - \left\| \sum_{i\in J_{1} \setminus K} \left( - x_i \right) \right\|  \right|  \\
 &\leq \left\| \sum_{i\in K } x_i - \sum_{i\in J_{1} \setminus K} \left( - x_i \right)  \right\| \\
 &=  \left\| \sum_{i\in K } x_i + \sum_{i\in J_{1} \setminus K}  x_i   \right\| \\
 &= \left\|  \sum_{i\in K \cup J_{1}} x_i \right\| \\
 & \leq 1 + \left\| s \right\|.
\end{align*}
Then, $$ \displaystyle \left\| \sum_{i\in K } x_i \right\| \leq  1 + \left\| s \right\| + 
\left\| \sum_{i\in J_{1} \setminus K} x_i \right\| \leq 1 + \left\| s \right\| + 
 \sum_{i\in J_{1} \setminus K} \left\|  x_i \right\| .$$
Hence,  
$$\displaystyle \left\| \sum_{i\in K } x_i \right\| \leq  1 + \left\| s \right\| + 
\sum_{i\in J_{1}} \left\|  x_i \right\|  $$
which completes the proof. 
\end{proof}

The following proposition establishes that norm summability forces a net to converge to zero in norm whenever its directed index set has no maximal element.
\begin{proposition} \label{p2}
Let $(x_\alpha)_{\alpha \in I}$ be a net in  $X$.
Suppose that the directed index set $I$ has no maximal element.\\
 If 
${\left( x_{\alpha}\right)_{\alpha \in I}}$ 
is a norm summable family, then $x_{\alpha}\xrightarrow{\|.\|} 0$.
\end{proposition}

\begin{proof}
Let $(x_\alpha)_{\alpha \in I}$ be a norm summable family in $X$ such that \\
$\displaystyle \sum_{\alpha \in I} {x_ \alpha}_{ / \|.\|} = s$, and let $\varepsilon$ be a strictly positive real number, then there exists a finite subset $J_{\varepsilon} \subseteq I$ such that for every finite subset 
$J \subseteq I$, $$J_\varepsilon \subseteq J \implies \left\| \displaystyle \sum_{i \in J} x_i - s \right\| \le \dfrac{\varepsilon}{2}.$$
On the other hand, since $I$ has no maximal element, there exists 
$\alpha_0 \in I$ such that $  \alpha_0  > i \quad \text{for every } i \in J_{\varepsilon}$.
 So, for every 
$\alpha \ge \alpha_0$, we have 
\begin{align*}
\| x_{\alpha}\| &= \left\| \displaystyle \sum_{i \in J_{\varepsilon}\cup \{\alpha\}} x_i - \displaystyle \sum_{i \in J_{\varepsilon}} x_i \right\| \\
  &\leq \left\| \displaystyle \sum_{i \in J_{\varepsilon}\cup \{\alpha\}} x_i - s \right\| + 
  \left\| \displaystyle \sum_{i \in J_{\varepsilon}} x_i - s \right\| \\
  &\leq \dfrac{\varepsilon}{2}+\dfrac{\varepsilon}{2}= \varepsilon
\end{align*}
Therefore, $x_{\alpha}\xrightarrow{\|.\|} 0$.
\end{proof}

\begin{corollary}
Let $(x_\alpha)_{\alpha \in I}$ be a norm summable net in a normed vector space $X$.
Then $x_{\alpha}\xrightarrow{\|.\|} 0$  if and only if one of the following conditions is satisfied:
\begin{enumerate}[label=(\arabic*)]
\item The index set $I$ has no maximal element.
\item The index set $I$ has a maximal element $\alpha_0$ and $x_{\alpha_0}=0$.
\end{enumerate}
\end{corollary}

\begin{lemma} \label{l1}
Let $X$ be a normed Riesz space and let ${(x_\alpha)}_{\alpha \in I}$ be a net in $X$.
If $x_\alpha \xrightarrow{o} a$ and $x_{\alpha}\xrightarrow{\|.\|} b$,
then $a = b$.
\end{lemma}

\begin{proof}
Let $x_\alpha \xrightarrow{o} a$, and $x_{\alpha}\xrightarrow{\|.\|} b$.\\
Then there exists ${(y_{\beta})}_{\beta \in J}$ in $X^{+}$ satisfying: $y_{\beta} \downarrow 0$ and for any 
$\beta \in J$, there exists ${\alpha}_{0} \in I$ such that 
${\left \lvert x_{\alpha} -a  \right \rvert} \leq y_{\beta}$  for all $\alpha \geq {\alpha}_{0}$.\\
Now fix $\alpha \geq {\alpha}_{0} $ in $I$. So 
$$\left| a-b \right| \leq \left| x_\alpha - b \right| + \left| x_\alpha -a \right|$$ which implies that 
$$\left| a-b \right| \leq \left| x_\alpha - b \right| + \left| y_\beta \right| \xrightarrow{o} \left| x_\alpha - b \right|$$ 
Therefore $\left| a-b \right| \leq \left| x_\alpha - b \right|$ for all $\alpha \geq {\alpha}_{0} $ in $I$, then
$$\left\| a-b \right\| \leq \left\| x_\alpha-b \right\| \xrightarrow{} 0$$
Hence $a=b$.
\end{proof}

The following proposition establishes the compatibility between semi-order summability and norm summability. It shows that whenever a family is summable in both senses, the corresponding sums necessarily coincide.
\begin{proposition} \label{p8}
Let $X$ be a normed Riesz space. If ${(x_i)}_{i \in I}$ is a family in $X$ that is both semi-order summable and norm summable, then 
$$\displaystyle \sum_{i \in I} {x_i}_{ / so}= \displaystyle \sum_{i \in I} {x_i}_{ / \|.\|}.$$
\end{proposition}

\begin{proof}
We assume that ${(x_i)}_{i\in I}$ is a semi-order summable and norm summable family in $X$,
$\displaystyle \sum_{i \in I} {x_i}_{ / so}=s$, and that 
$\displaystyle \sum_{i \in I} {x_i}_{ / \|.\|}=s'$.\\
Let $\Gamma=\biggl\{ \displaystyle \beta \mid \beta \subseteq I \text{ with } \beta \text{ is finite } \biggl\}$ and $X_{\beta} = \displaystyle \sum_{i \in \beta} {x_i} ,\; \beta\in\Gamma$. 
So, 
$$X_{\beta} \xrightarrow{o} s  \quad \text{and} \quad X_{\beta} \xrightarrow{\|.\|} s'$$
By applying Lemma~\ref{l1}, we have $s = s'$.

\end{proof}

Unlike in the case of summability in normed vector spaces, one can find order summable families with uncountable support in a normed Riesz space. Consequently, there are families that are order summable but not norm summable. The following example illustrates this.

\begin{example}
Let $I=\left]0,1\right]$, $E=\displaystyle \prod_{i\in I} \mathbb{R}$ and 
$$
E_{\infty} = 
\Bigl\{\, x = (x_i)_{i \in I} \in E \ \Big|\ 
\exists M \in \mathbb{R}_+ \text{ such that } |x_i| \le M \text{ for all } i \in I \,\Bigr\}
$$
We equip $E_{\infty}$ with its usual order and the norm 
$\|.\|_{\infty}: E_{\infty} \longrightarrow \mathbb{R}_+$ defined by 
$\|x\|_{\infty}=\displaystyle \sup_{i \in I} |x_i|$ for all $x = (x_i)_{i \in I} \in E_{\infty}$.\\
Now, for each $j \in I$, we define the map 
$ p_j \colon E_{\infty} \to E_{\infty} $
by
$ p_j\bigl( (a_i)_{i \in I} \bigr) = b_j, $
where $b_j = (b_j^i)_{i \in I}$ is defined as
$$
b_j^i =
\begin{cases}
a_j, & \text{if } i = j,\\
0, & \text{if } i \neq j.
\end{cases}
$$ 

Then, the family $(p_i(x))_{i \in I}$ is order summable, and $\displaystyle \sum_{i \in I} {p_i(x)}_{ / o} = x$ for each $x \in E.$ However, the family $(p_i(x))_{i \in I}$ is not norm summable if \\ $x ={\left( i\right)}_{i \in I}$, since $I$ is uncountable and  $p_i(x) > 0$ for each $i \in I$. 
\end{example}

The next proposition states that any family that is semi-order summable is also norm summable, provided that the underlying normed Riesz space has an order-continuous norm.
\begin{proposition} \label{p9}
Let $X$ be a normed Riesz spacee with order-continuous norm. Then every semi-order summable family in $X$ is also norm summable.
\end{proposition}
 \begin{proof}
 We assume that  $(x_i)_{i \in I}$ be a semi-order summable family in $X$, and that
 $\displaystyle \sum_{i \in I} {x_i}_{ / so}=s$ \\ 
 Let $\Gamma=\biggl\{ \displaystyle \beta \mid \beta \subseteq I \text{ with } \beta \text{ is finite } \biggl\}$, and let $X_{\beta} = \displaystyle \sum_{i \in \beta} {x_i}, \; \beta \in \Gamma$.
So,  $ X_{\beta} \xrightarrow{o} s$. 
Since the norm of $X$ is order-continuous, it follows that 
$ \left( X_{\beta}-s  \right) \xrightarrow{\|.\|} 0$. 
Consequently, $X_{\beta} \xrightarrow{\|.\|} s$, and therefore 
 $(x_i)_{i \in I}$ is norm summable.
 \end{proof}
 
 \begin{corollary}
Let $X$ be a normed Riesz spacee with order-continuous norm. Then every order summable family in $X$ is also norm summable.
\end{corollary}
\begin{corollary}
Let $(x_i)_{i \in I}$ be a semi-order summable family in a normed Riesz space $X$ with order-continuous norm. Then the set 
 $\{i \in I \mid x_i \neq 0\}$ is countable.
\end{corollary}

\begin{proof}
Let $(x_i)_{i \in I}$ be a semi-order summable family in $X$. Since $X$ is a normed Riesz space, it follows from Proposition~\ref{p9} that $(x_i)_{i \in I}$ is a norm summable family. Hence, by Remark~\ref{r2}, the set 
 $\{i \in I \mid x_i \neq 0\}$ is countable.
\end{proof}

The following theorem establishes that the order continuity of the norm in a Banach lattice can be characterized by the interplay between order summability and norm summability on the one hand, and between semi-order summability and norm summability on the other.
\begin{theorem} \label{t2}
Let $X$ be a Banach lattice. Then the following assertions are equivalent:
\begin{enumerate}[label=(\arabic*)]
    \item $X$ has order-continuous norm.
    \item $X$ is Dedekind complete and every semi-order summable family in $X$ is also norm summable.
    \item $X$ is Dedekind complete and every order summable family in $X$ is also norm summable.
    \item $X$ is Dedekind $\sigma$-complete and every semi-order summable sequence in $X$ is also norm summable.
    \item $X$ is Dedekind $\sigma$-complete and every order summable sequence in $X$ is also norm summable.
\end{enumerate}
\end{theorem}

\begin{proof}
$\text{(1)} \Rightarrow \text{(2)}$ by \cite[Corollary~4.10]{cd1} and Proposition~\ref{p9}.\\
$\text{(2)} \Leftrightarrow \text{(3)}$ by Theorem~\ref{t6}.\\
$\text{(4)} \Leftrightarrow \text{(5)}$ by Theorem~\ref{pp8}.\\
$\text{(3)} \Rightarrow \text{(5)}$ is trivial.\\
$\text{(5)} \Rightarrow \text{(1)}$ 
Let $(x_n)_{n \in \mathbb{N}}$ be an order bounded disjoint sequence in $X$. 
By Corollary~\ref{cr1}, $(x_n)_{n \in \mathbb{N}}$ is order summable. 
Hence, it is norm summable. 
By Proposition~\ref{p2}, we obtain ${x_n} \xrightarrow{\|.\|} 0$.
Therefore, by \cite[Theorem~4.14]{cd1}, $X$ has order-continuous norm.
\end{proof}

Let $(x_i)_{i \in I}$ be a family in a normed Riesz space $X$. 
If the families $({x}_{i}^{+})_{i \in I}$ and $({x}_{i}^{-})_{i \in I}$ 
are norm summable in $X$, then $(x_i)_{i \in I}$ is also 
norm-summable in $X$. However, in general, the converse fails: 
it may happen that $(x_i)_{i \in I}$ is norm summable whereas 
neither $({x}_{i}^{+})_{i \in I}$ nor $({x}_{i}^{-})_{i \in I}$ is 
norm summable, as shown by the following example.

\begin{example}
Let 
$X = c_{00}$, endowed with its usual order and the norm defined, for
$x=(x_{1},x_{2},\dots)\in X,$
by
$\|x\|=\sup\{|x_i|: i\in \mathbb{N}\}.$
Let $({x}_{n})_{n \in \mathbb{N}^*}$ be the family in $X$ defined by
${x}_{n}=\dfrac{1}{n}e_{n}-\dfrac{1}{n+1}e_{n+1},$
where $(e_{n})_{n \in \mathbb{N}^*}$ denotes the canonical basis vectors.
Then $({x}_{n})_{n \in \mathbb{N}^{*}}$ is norm summable in $X$ and 
$\displaystyle \sum_{n \in \mathbb{N}^{*}} {x_{n}}_{ / \|.\|}=e_{1}$. However,  
$({x}_{n}^{+})_{n\in\mathbb{N}^*}
\quad \text{and} \quad
({x}_{n}^{-})_{n\in\mathbb{N}^*}$
are not norm summable in $X$.
\end{example}

The following theorem establishes that, in a KB-space, the norm summability of a disjoint family is completely determined by the norm summability of its positive and negative parts, or equivalently, by the norm summability of the family of its moduli.
\begin{theorem} \label{tkb}
Let $X$ be a KB-space, and let $(x_i)_{i\in I}$ be a disjoint family in $X$. Then the following statements are equivalent:
\begin{enumerate}
\item There exists a constant $M \in \mathbb{R}_+$ such that:
$$\left\Vert{} \sum_{i \in \beta} x_i \right\Vert{} \leq M \quad \text{for every finite subset } \beta \subseteq I$$
\item $(x_i)_{i\in I}$ is norm summable in $X$.
\item the families $(x_i^{+})_{i\in I}$ and $(x_i^{-})_{i\in I}$ are norm summable in $X$.
\item the family $(|x_i|)_{i\in I}$ is norm summable in $X$.
\end{enumerate}
\end{theorem}

\begin{proof}
Let 
$\Lambda=\biggl\{ \displaystyle \beta \mid \beta \subseteq I \text{ with } \beta \text{ is finite } \biggl\}$,
$X_{\beta}= \displaystyle \sum_{i \in \beta } {x_i}^{+}$ and  
$Y_{\beta}= \displaystyle \sum_{i \in \beta } {x_i}^{-}$. 
The implications $(3)\Rightarrow(2)$ and $(3)\Rightarrow(4)$ are immediate by Propositin \ref{p1}.
The implication $(2)\Rightarrow(1)$ by Proposition~\ref{l5}\\
$(1)\Rightarrow(3)$ 
Showing that $({x_i}^{+})_{i\in I}$ is norm summable in $X$ is equivalent to showing that the net ${(X_{\beta})}_{\beta \in \Lambda}$ is norm convergent in $X$.\\
First, it is clear that the net ${(X_{\beta})}_{\beta \in \Lambda}$ is increasing. Since $X$ is a KB-space, it is enough to show that 
${(X_{\beta})}_{\beta \in \Lambda}$ is norm bounded.
Since $(x_i){i\in I}$ is a disjoint family, the family $(x_i^+){i\in I}$ is also disjoint. Hence, for every finite subset $\beta\subseteq I$,
$$X_{\beta}=\sum_{i\in\beta}x_i^+
=\left(\sum_{i\in\beta}x_i\right)^+
\leq \left|\sum_{i\in\beta}x_i\right|.$$
So 
$$\displaystyle \| X_{\beta} \| \leq  \left\| \sum_{i\in\beta}x_i \right\| \leq M $$

Similarly, one shows that $(x_i^{-})_{i\in I}$ is norm summable in $X$.\\
$(4)\Rightarrow(3)$  
Assume that the family $(|x_i|)_{i\in I}$ is norm summable. Then, by Proposition~\ref{l5}, there exists 
$M \in \mathbb{R}^{+}$ such that
$\displaystyle \left\|\sum_{i\in \beta} |x_i|\right\|\leq M$
for each $\beta \in \Lambda$.\\
Since
$$X_{\beta} \leq \sum_{i\in \beta} |x_i| \quad \text{and} \quad Y_{\beta} \leq \sum_{i\in \beta} |x_i|,$$
it follows that
$$\|X_{\beta} \| \le M \quad \text{and} \quad \| Y_{\beta} \| \le M.$$
It is clear that ${(X_{\beta})}_{\beta \in \Lambda}$ and ${(Y_{\beta})}_{\beta \in \Lambda}$ are increasing nets. Since $X$ is a KB-space, it follows that both ${(X_{\beta})}_{\beta \in \Lambda}$ and ${(Y_{\beta})}_{\beta \in \Lambda}$ are norm convergent. Thus, $(x_i^{+})_{i\in I}$ and $(x_i^{-})_{i\in I}$ are norm summable in $X$.
\end{proof}
\textbf{Question:} Under the same conditions on the space as in the previous theorem, do the equivalences remain valid for nets that are not disjoint?
\section{Conclusion}
In this paper, we introduced the notion of semi-order summability for arbitrary families in Riesz spaces and investigated its relationship with order summability and norm summability. We established the fundamental properties of semi-order summable families and demonstrated that semi-order summability and order summability are equivalent under certain assumptions.. We also obtained new characterizations of laterally complete Riesz spaces and Banach lattices with order-continuous norm in terms of summability, and established a characterization of norm summability of disjoint families in KB-spaces. 
These results provide a unified framework for the study of summability of arbitrary families in Riesz spaces and Banach lattices and suggest several directions for future research, including the study of semi-order summability in more general ordered structures and of operators preserving these notions of summability.

\section*{Declarations}

\subsection*{Competing Interests}
The authors declare that they have no competing interests.

\subsection*{Funding Information}
Not Applicable.

\subsection*{Author Contribution}
All authors contributed equally to the conception, methodology, investigation, 
formal analysis, writing, and revision of the manuscript. All authors read 
and approved the final manuscript.

\subsection*{Data Availability Statement}
Not Applicable.


\end{document}